\input{epsf}

\documentclass[a4paper,12pt]{article}

\usepackage[latin1]{inputenc} 
\usepackage[T1]{fontenc}      
\usepackage[french]{babel}
\usepackage{fullpage}
\usepackage{amssymb,amsfonts}
\usepackage{graphics}
\usepackage{caption}
\usepackage{amsthm,amsmath}
\usepackage{hyperref}

\def\S{\mathbb{S}}

\def\R{\mathbb{R}}
\def\C{\mathbb{C}}

\def\Z{\mathbb{Z}}
\def\OO{\mathcal{O}}
\def\cI{\mathcal{I}}
\def\II{\mathrm{I\hskip-1pt I}}
\def\im{\mathrm{im}\,}
\def\ker{\mathrm{ker}\,}
\def\rg{\mathrm{rg}\,}

\def\epsilon{\varepsilon}
\def\phi{\varphi}

\newtheorem{thm}{Théorème}
\newtheorem{prop}[thm]{Proposition}
\newtheorem{lemm}[thm]{Lemme}
\newtheorem{cor}[thm]{Corollaire}

\theoremstyle{definition}
\newtheorem{exemple}{Exemple}
\newtheorem*{defn}{Définition}
\newtheorem*{remq}{Remarque}

\def\cqfd{\hfill $\Box$\break\null}
\def\dem{\noindent{\it Démonstration.\  }}

\begin{document}

\title{Les douze surfaces de Darboux et la trialité}
\author{Bruno Sévennec\thanks{C.N.R.S., U.M.P.A., \'{E}cole Normale Supérieure de Lyon, 46 All\'ee d'Italie, 69364 Lyon cedex 07}
}

\maketitle

\begin{abstract}
 On donne dans cet article une interprétation géométrique des "douze surfaces de Darboux", qui apparaissent en appliquant de façon répétée une transformation simple à une déformation isométrique infinitésimale d'une surface dans l'espace euclidien de dimension trois. Cette interprétation est une version différentielle de la trialité, concernant les immersions totalement isotropes de surfaces dans la quadrique projective réelle de dimension 6 définie par une forme quadratique de signature neutre (4,4). 
\end{abstract}


\section{Introduction}

Dans le volume 4 du traité de Gaston Darboux sur les surfaces \cite{Darboux}, 
consacré aux déformations isométriques des surfaces dans l'espace euclidien $\R^3$,
on trouve au chapitre 3, intitulé ``les douze surfaces'', une construction que
l'on peut présenter comme suit (voir aussi \cite{Sabitov}).

Soit $S$ une surface, $f,g$ deux applications (lisses)
de $S$ vers $\R^3$ avec $f$ une immersion ($df$ est partout injective), telles que 
\begin{equation}
	 df\cdot dg=0, 
	\label{eq:1}
\end{equation}
i.e. qu'en tout point $p$ de $S$ on ait pour tout $v\in T_p S$ 
$$df(p)v \cdot dg(p)v = 0,$$
(on note $\cdot$ le produit scalaire euclidien dans $\R^3$).
Cette condition est équivalente au fait que la métrique riemannienne $q_\epsilon$
induite sur $S$ par $f_\epsilon=f+\epsilon g : S \to \R^3_\mathrm{eucl}$ vérifie
$$q_\epsilon=q_0 +O(\epsilon^2).$$
Autrement dit (\ref{eq:1}) exprime que $g$, vu comme champ de vecteurs le long de $f$,
est une déformation isométrique infinitésimale (d'ordre 1) de $f:S\to\R^3$, sujet
de nombreuses investigations (notamment en ingénierie) au temps de Darboux comme par la suite.

Du fait que $f$ est une immersion, il existe une unique
application (lisse) $h:S\to\R^3$ telle que 
\begin{equation}\label{eq:2}
dg(p) = h(p) \times df(p),\ \ p\in S
\end{equation}
où $\times$ désigne le
produit vectoriel dans $\R^3$. On peut concevoir $h$ comme le champ des ``rotations
infinitésimales'' induites par la déformation $g$ de $f$.
Darboux remarque alors que $\tilde{g}=g-h\times f$ vérifie $$d \tilde{g} = f\times dh.$$

Si de plus $g$ est une immersion, la relation (\ref{eq:1}) étant symétrique en $f,g$ 
on peut définir $h^*: S \to \R^3$ par $$df = h^* \times dg.$$  

Il est clair que ces formules définissent, sur les triplets $(f,g,h)$ vérifiant (\ref{eq:2})
ainsi que des conditions de non-dégénérescence convenables, deux  involutions
$$
\begin{array}{cccc}
	 A:& (f,g,h)& \mapsto &(h,\tilde{g},f) \\
   D:& (f,g,h)& \mapsto &(g,f,h^*) \;\;.
\end{array}
$$

Le résultat surprenant obtenu par Darboux peut maintenant s'énoncer
\begin{thm}[Darboux] \label{th:ordre6}
La composition $(D\circ A)^6$ est l'identité.
\end{thm}
Les ``douze surfaces'' de Darboux sont les premières composantes des douze triplets
obtenus par application des éléments du groupe diédral engendré par $A$ et $D$ 
au triplet initial $(f,g,h)$.

On propose ici une formulation géométrique de ce théorème, sous forme
d'une ``trialité différentielle'' pour les surfaces totalement isotropes
dans la quadrique projective réelle $Q \subset \R P^7 = P(\R^8)$ 
définie par l'annulation d'une forme quadratique de signature $(4,4)$ sur $\R^8$.

Le lien entre la propriété (\ref{eq:1}) et la quadrique $Q$ dans $\R P^7$ vient de ce que,
considérée comme propriété de l'application
$$(f,g):S \to \R^3\times\R^3$$
elle ne dépend visiblement que de la structure (pseudo-) conforme définie par la métrique pseudo-riemannienne
plate $dx\cdot dy$ de signature $(3,3)$ sur $\R^6=\R_x^3\times\R_y^3$. 

Or, par une construction bien connue (rappelée dans la section \ref{secff}) généralisant la projection stéréographique usuelle,
la quadrique $Q\subset \R P^7$ est une compactification (pseudo-) conforme de $(\R^6,dx\cdot dy)$.
Plus précisément, si on munit $\R^8=\R^6\times\R^2=\R_x^3\times\R_s\times\R_y^3\times\R_t$ 
de la forme quadratique $$q(x,s,y,t)=x\cdot y +st,$$  
on identifie $\R^3\times\R^3$ à l'ouvert affine défini par $s=1$ dans la quadrique projective $Q$
d'équation $q=0$ dans $\R P^7$ au moyen de l'application 
$$(x,y)\mapsto (x:1:y:-x\cdot y)$$
(rappelons que $x\cdot y$ désigne le produit scalaire dans $\R^3$).
\footnote{
Le cas usuel correspond à la quadrique projective de signature $(n+1,1)$ dans $\R P^{n+1}=P(\R^{n+2})$,
i.e. $\S^n$, qui est une compactification conforme de $\R^n_\mathrm{eucl}$}

Un couple $(f,g)$ vérifiant (\ref{eq:1}) devient alors une immersion $\phi : S \to Q$ totalement isotrope, 
i.e. dont le plan projectif tangent en tout point est entièrement contenu dans la quadrique $Q$, et 
on voit facilement qu'il est alors inclus dans exactement {\em deux} espaces projectifs maximaux
(de dimension $3$) contenus dans $Q$.
 
Or il se trouve que les espaces projectifs de dimension $3$ contenus dans $Q$ forment deux familles $Q_\pm$, elles-mêmes
isomorphes (algébriquement) à la quadrique $Q$, et on associe ainsi à $\phi$ {\em deux} applications $\phi_\pm:S\to Q_\pm$, 
elles aussi lisses et totalement isotropes, mais qui peuvent dégénérer, i.e. être de rang $<2$.

Dans cette situation, la trialité projective est le fait 
que les relations d'incidence naturelles entre éléments de $Q,Q_+,Q_-$ définissent sur leur union disjointe
une structure hautement symétrique, dont les automorphismes permutent arbitrairement les trois quadriques. 
Ce fait est classique, voir \cite{vdBS} et ses références, dont \cite{Cartan}.

La ``trialité différentielle'' mise en évidence ici affirme que $\phi: S \to Q$ et
ses deux ``applications de Gauss'' $\phi_\pm : S \to Q_\pm$ jouent aussi des rôles
parfaitement symétriques (au moins lorsque $\phi$ est générique), 
et on montrera que ceci permet de retrouver
le théorème \ref{th:ordre6} de Darboux.

Dans la suite, on commence par rappeler sections \ref{secff} à \ref{dii} 
les faits ``classiques''\footnote{mais méritant peut-être d'être rappelés}
ci-dessus, ainsi que quelques autres qui leurs sont reliés, 
comme la notion de seconde forme fondamentale projective, 
la trialité sur la quadrique projective "neutre" de dimension $6$ et le théorème de Darboux-Sauer 
sur l'invariance projective de la rigidité infinitésimale. 

Puis on étudie section \ref{triplets} les triplets de Darboux et
leurs involutions $A$, $D$, préparant le calcul section \ref{calcul}
des deux ``applications de Gauss'' $\phi_\pm$ de l'immersion totalement isotrope $\phi$ 
associée à une déformation isométrique infinitésimale $(f,g)$. Génériquement, elles s'expriment
en termes de composantes de triplets de Darboux transformés du triplet initial $(f,g,h)$.

Après avoir introduit section \ref{octo} le modèle des matrices de Zorn pour les octonions déployés, on met en
évidence la trialité différentielle section \ref{trialdiff}, puis on identifie section \ref{darbouxZorn} 
les tranformations de Darboux en termes de cette trialité, démontrant au passage le théorème \ref{th:ordre6}.
On caractérise ensuite section \ref{degen} les cas dégénérés où l'une des applications $\phi_\pm$ est partout de rang $<2$.
Les trois ``secondes formes fondamentales'' associées 
à une immersion totalement isotrope sont définies et étudiées section \ref{secff3}.

Enfin on considère section \ref{varincid} une variété d'incidence de dimension $11$, 
munie d'un champ de plans tangents de dimension $6$ complètement non-intégrable,
dont les immersions totalement isotropes sont des surfaces intégrales.

Remarquons pour finir que ces constructions pourraient permettre d'aborder le
problème toujours ouvert (voir \cite{Yau,Ghomi}) de l'existence de surfaces lisses compactes 
sans bord (nécessairement non convexes) dans $\R^3$ admettant des 
déformations isométriques lisses (finies) non triviales. 
Rappelons qu'il existe des exemples polyédraux \cite{Connelly,Kuiper}. 
On peut ainsi se demander s'il existe, dans la catégorie \emph{analytique réelle} $C^\omega$, 
des surfaces compactes sans bord admettant des déformations isométriques infinitésimales non triviales.
Trotsenko montre dans \cite{Trotsenko} qu'il existe des surfaces de révolution 
difféomorphes à $S^2$ possédant des déformations isométriques
infinitésimales non triviales, et même un espace 
de dimension arbitrairement grande de telles déformations, 
mais comme le fait remarquer Spivak dans \cite{Spivak},
ces déformations ne sont pas lisses aux pôles~: elles ont un ordre
de différentiabilité fini.


\section{La seconde forme fondamentale d'une hypersurface dans l'espace projectif}\label{secff}

Soit $M\subset \R P^n$ une hypersurface de classe $C^2$, $p$ un point de $M$, et $f=0$ une équation locale $C^2$ de $M$ au voisinage de $p$.
Autrement dit, $M$ coïncide avec $f^{-1}(0)$ au voisinage de $p$, et $df(p)\neq 0$.

Il existe un unique hyperplan projectif $H$ tangent à $M$ en $p$, i.e. tel que $T_pH=T_pM$ (il est déterminé par $\ker df(p)$), et la restriction $\overline{f}=f|_H$ à $H$ possède un point critique en $p$. Il est alors classique (et facile à prouver) que $\overline f$ admet en $p$ une différentielle seconde intrinsèque, forme bilinéaire symétrique $D^2\overline{f}(p)\in S^2T^*_pH=S^2T^*_pM$.

Puisque $df(p)\neq 0$, il existe une unique application bilinéaire symétrique $$\II_p : S^2 T_pM \to N_pM := T_p(\R P^n) /T_pM$$
telle que sa composée avec $$df(p):N_pM \tilde{\to} \R$$
soit $-D^2\overline{f}(p)$. Si on change d'équation locale de $M$ en $p$, $df(p)$ et $D^2\overline{f}(p)$ sont multipliés par un même facteur non nul, et $\II_p$ ne dépend d'aucun choix. C'est la seconde forme fondamentale de $M$ en $p$. Elle définit un sous-espace vectoriel de dimension au plus $1$ dans $S^2T^*_pM$ (via l'ensemble des isomorphismes $N_pM\simeq \R$), que l'on appellera structure pseudo-conforme projective de $M$ en $p$.

\begin{exemple}
Soit $q$ une forme quadratique sur $\R^n$, et $M\subset \R^{n+1}\subset\R P^{n+1}$ le graphe $\{t=q(x)\}$. 
Alors la structure pseudo-conforme projective de $M$ s'identifie via la projection $M\tilde{\to}\R^n$ à celle associée au champ constant de formes quadratiques $$x\in\R^n \mapsto \II_x=2q\in S^2(\R^n)^*\simeq S^2T^*_x\R^n.$$
Cela se vérifie sans difficulté, avec $f(x,t)=q(x)-t$ et $N_pM$ identifié à l'axe vertical, 
mais peut aussi se déduire de l'existence d'une action affine transitive sur $M\subset\R^{n+1}$ 
qui relève celle des translations sur $\R^n$~: $$(x,t)\mapsto (x+v,t+dq(x)v+q(v)),\ (x,t)\in\R^{n+1},v\in\R^n.$$
Cette action est bien sûr restriction d'une action projective, et ne reste plus que le calcul trivial en $(x,y)=(0,0)$.
\end{exemple} 

\begin{exemple}
Dans l'exemple précédent, si on suppose $q$ non dégénérée, $M$ est contenue dans la quadrique projective lisse $Q\subset\R P^{n+1}$
d'équation $q(x)-st=0$ en coordonnées homogènes $(s:x:t)$. Plus précisément, $M=Q\setminus H$, 
où $H$ est l'hyperplan à l'infini $s=0$, qui est tangent à $Q$ au point à l'infini sur l'axe vertical $p=(0:0:1)$.
Prenant $q(x)=|x|^2$ euclidienne, et $H'$ l'hyperplan $s+t=0$, on voit que $Q$ s'identifie à une sphère dans l'espace affine
$\R P^{n+1}\setminus H'$, et que la projection stéréographique $\pi:Q\setminus\{p\}\to \R^n$ depuis $p$ s'identifie à la projection verticale de l'exemple précédent. On retrouve ainsi le fait que $\pi$ est conforme, vu la définition projective de (la classe conforme de) la seconde forme fondamentale. 
\end{exemple}

Ces exemples montrent que tout espace affine réel muni d'une métrique pseudo-riemannienne plate de signature $(p,q)$ admet une compactification (pseudo-) conforme naturelle, qui est une quadrique projective lisse de signature $(p+1,q+1)$ dans un espace projectif de dimension $p+q+1$.
Le ``lieu à l'infini'' qui est ajouté est le cône sur une quadrique projective de signature $(p,q)$ (un simple point si $p$ ou $q$ est nul).


\section{Trialité}\label{trial}
Parmi les algèbres de Lie simples complexes, celle de type $D_4$, autrement dit $\mathfrak{so}(8,\C)$ est exceptionnelle en ce qu'elle admet un automorphisme extérieur d'ordre $3$, qui permute cycliquement les trois branches (de longueur 1) de son diagramme de Dynkin en forme de Y. 
Ce fait peut paraître bien abstrait, mais il 
s'incarne géométriquement de manière très concrète dans la quadrique réelle $Q$ de signature $(4,4)$  dans $\R P^7$, espace homogène du groupe de Lie réel $PO(4,4)$ de type $D_4$.

Plus précisément, considérons l'espace vectoriel réel $\R^{4,4}=\R^4\oplus\R^4$ muni de la forme quadratique $q(u,v)=|u|^2-|v|^2$ de signature $(4,4)$. Alors $Q$ est l'espace des droites vectorielles isotropes pour $q$.

Un sous-espace isotrope maximal (de dimension $4$) de  $\R^{4,4}$ n'est autre que le graphe d'une isométrie euclidienne $\R^4\to\R^4$, et la variété de ces sous-espaces a donc deux composantes connexes, identifiées à celles $O^\pm(4)$ du groupe réel $O(4)$.

Or $SO(4)=O^+(4)$ est classiquement isomorphe  à $S^3\times S^3/{\pm (1,1)}$, en voyant $S^3$ comme le groupe multiplicatif des quaternions de norme $1$,
agissant par multiplications à gauche et à droite sur l'algèbre des quaternions identifiée à $\R^4$.

Mais $S^3\times S^3/{\pm (1,1)}$ est précisément une copie de la quadrique projective $Q$ d'équation homogène $|u|^2-|v|^2=0$. Ainsi les deux familles de $\R P^3$ contenues dans $Q$ sont deux copies $Q_\pm$ de la quadrique $Q$.
  
Il existe entre deux quelconques de ces quadriques une relation d'incidence. 
Entre $Q$ et $Q_\pm$, c'est simplement la relation d'appartenance, et entre $Q_+$ et $Q_-$, 
c'est la condition d'avoir une intersection maximale, autrement dit un plan projectif (contenu dans $Q$).
 
En termes d'isométries euclidiennes $g_\pm\in O^\pm(4)$, cela signifie que $g_+$ et $g_-$ coïncident
sur un sous-espace de dimension $3$, autrement dit diffèrent par une réflexion hyperplane 
(on voit ainsi que les éléments de $Q_-$ incidents à un élément fixé de $Q_+$ forment un $\R P^3$, et vice-versa).

De même, on voit que tout plan projectif contenu dans $Q$ est inclus dans exactement un élément de $Q_+$ et un élément de $Q_-$, éléments dont il est l'intersection.
Le fait qu'un sous-espace vectoriel $V$ de $\R^{4,4}$ totalement isotrope de dimension $3$ est contenu dans exactement deux sous-espaces
totalement isotropes maximaux peut aussi se voir simplement en remarquant que $V^\bot/V\simeq \R^{1,1}$ a exactement deux droites isotropes.

Le résultat suivant est plus profond
\begin{thm} Le groupe des bijections de l'union disjointe $Q\cup Q_+\cup Q_-$ qui respectent les relations d'incidence agit sur les composantes en les permutant arbitrairement.
\end{thm}
C'est l'origine du terme trialité, en référence au fait analogue pour l'union d'un espace projectif et de son dual, à savoir la dualité (projective).

La preuve la plus simple de ce théorème (voir \cite{vdBS}) consiste à introduire sur $\R^{4,4}$ une structure d'algèbre (non associative, à unité), celle des octonions déployés, pour laquelle la forme quadratique $q$ est multiplicative. On montre alors que $Q_\pm$ sont paramétrés projectivement par les octonions de norme nulle $a$, via leurs annulateurs gauche et droit, et que les relations d'incidence entre $a,a_+,a_-$ sont simplement
$aa_-=0$, $a_-a_+=0$, $a_+a=0$. Le théorème en résulte aussitôt.

On utilisera par la suite un modèle particulier des octonions déployés, bien adapté aux calculs, celui des matrices de Zorn (voir section \ref{octo}).


\section{Déformations isométriques infinitésimales}\label{dii}

Soit $S$ une surface et $f:S\to \R^3$ une immersion dans l'espace euclidien, toutes deux lisses. 
Une déformation infinitésimale de $S$ est un champ de vecteurs
$g:S\to\R^3$ le long de $S$, et elle est isométrique si la dérivée en $t=0$ de la métrique riemannienne 
$|df_t|^2=|df|^2+2t\,df\cdot dg+t^2|dg|^2$ induite par $f_t=f+t\, g$ est nulle, autrement dit si $$df\cdot dg=0.$$
Un déplacement infinitésimal $g=a\times f +b$ (où $a,b\in \R^3$ et $\times$ désigne le produit vectoriel) est une déformation infinitésimale triviale.

Considérant le couple $\phi=(f,g)$ comme une immersion de $S$ dans 
l'espace pseudo-riemannien $\R^3\oplus\R^3$ muni de la 
métrique plate $dx\cdot dy$, cette condition dit simplement que 
$\phi$ est une immersion totalement isotrope~: l'image de $d\phi(p)$ 
est totalement isotrope pour tout point $p$ de $S$.

Cette condition ne dépendant évidemment que de la structure (pseudo-)conforme  
définie par $dx\cdot dy$, il est naturel d'après la section \ref{secff}
de considérer $\phi$ comme une immersion dans la quadrique projective réelle
de signature $(4,4)$ 
$$Q =\{ (x:s:y:t)\in\R P^7 \mid x\cdot y + st=0 \},$$
donnée par $$\phi=(f:1:g:-f\cdot g).$$
On notera $q$ la forme quadratique $q(x,s,y,t)=x\cdot y + st$, et 
$\R^{4,4}$ l'espace vectoriel $\R^4\oplus\R^4$ muni de cette forme. 

La condition d'isométrie infinitésimale sur $(f,g)$ est équivalente à la totale isotropie
de l'immersion $\phi:S\to Q$ ci-dessus pour la seconde forme fondamentale projective de $Q$,
et revient à demander que pour tout $p\in S$ le plan projectif $\tau_p\phi$ tangent à $\phi$ en $p$ 
soit contenu dans $Q$.

On sait alors (voir section \ref{trial}) que $\tau_p\phi$ est contenu dans exactement 
deux espaces projectifs de dimension $3$ contenus dans $Q$, soit $\phi_\pm(p)$,
d'où deux applications $$\phi_\pm : S \to Q_\pm.$$

Pour fixer les notations on désignera par $Q_+$ la composante contenant 
les graphes $(a,b)_+$ des applications linéaires (antisymétriques si on identifie les $\R^4$ source et but)
$$(x,s)\mapsto (y,t)=(a\times x + bs,-b\cdot x)$$
$a,b\in\R^3$, et par $Q_-$  
celle contenant les ``graphes tordus'' $(a,b)_-$ des applications
$$(x,t)\mapsto (y,s)=(a\times x + bt,-b\cdot x)$$
$a,b\in\R^3$. Noter que ces graphes décrivent des ouverts affines\footnote{Au sens de la géométrie algébrique : ils sont complémentaires de l'intersection avec un hyperplan dans un plongement projectif, ici le plongement de Plücker 
$Q_\pm\to \mathrm{Gr}_4(\R^8)\to P(\wedge^4\R^8)$.} 
(denses) de $Q_+$ et $Q_-$. 

Lorsque $\phi:S\to Q$ provient de $(f,g):S\to \R^{3,3}$ avec $g=a\times f+b$ \ un déplacement infinitésimal,
$\phi_+$ est constante, égale à $(a,b)_+$.

Un autre cas de dégénérescence se rencontre si $f:S\to\R^3$ prend ses valeurs dans un plan affine $P$ 
évitant l'origine et $g$ est un champ de vecteurs lisse quelconque normal à ce plan, de sorte qu'on a bien $df\cdot dg=0$. 
Alors on vérifie aisémént que $\phi_-$ est constante, égale à $(0,b)_-$ si $P$ a pour équation $b\cdot x=-1$.
On pourrait encore ajouter à $g$ une rotation infinitésimale $a\times f$ pour obtenir $\phi_-=(a,b)_-$.

Pour mieux visualiser la situation, remarquons que la condition $df\cdot dg=0$ 
équivaut au fait que les deux immersions $f\pm \epsilon g:S\to\R^3$ sont isométriques, 
pour un $\epsilon>0$ assez petit, puisque 
$$|df\pm \epsilon dg|^2=|df|^2+\epsilon^2|dg|^2\pm 2\epsilon df\cdot dg.$$

Lorsque $g$ est un déplacement infinitésimal, $f\pm \epsilon g$ diffèrent 
par une isométrie affine directe de $\R^3$, et $\phi_+$ est constante. 
Par contre dans le second cas, $f\pm\epsilon g$ diffèrent par 
la réflexion orthogonale par rapport à $P$, et c'est $\phi_-$ qui est constante, 
mais pas $\phi_+$ en général.

\begin{defn}
On dit qu'une immersion $f:S\to\R^3$ est infinitésimalement rigide si toute déformation isométrique infinitésimale $g$ est un déplacement infinitésimal $g=a\times f + b$, $a,b\in\R^3$. 
\end{defn}

Pour passer commodément du cadre projectif au cadre affine, remarquons
que la quadrique $Q$ contient deux espaces projectifs de dimension $3$ disjoints
$P_1=P(\R^4\oplus 0)$ et $P_2=P(0\oplus\R^4)$ et qu'on dispose de deux projections
$\pi_1:Q\setminus P_2\to P_1$, $\pi_2:Q\setminus P_1\to P_2$. En particulier si $\phi=(f:1:g:-f\cdot g)$,
$\pi_1\circ\phi=(f:1)$ s'identifie à $f:S\to \R^3\subset P_1\simeq \R P^3$. De plus, la projection $\pi_1$
possède une structure naturelle de fibré vectoriel qui l'identifie au fibré $\mathrm{Hom}(L^\bot,L)$, 
où $L\subset P_1\times\R^4$ est le fibré en droites tautologique, de sorte que $\pi_1$ s'identifie à la projection
du fibré (co-)tangent à $P_1\simeq\R P^3$.

Alors une immersion $f:S\to\R^3$ est infinitésimalement rigide si et seulement si toute immersion 
totalement isotrope $\phi:S\to Q$ relevant $f$, i.e. telle que $\pi_1\circ\phi=f:S\to\R^3\subset\R P^3$,
est telle que $\phi_+$ est constante.

Un résultat classique est que c'est le cas lorsque $f$ est un plongement de la sphère à courbure de Gauss strictement positive. 
Au contraire, ce n'est pas le cas dès que $f$ envoie un ouvert dans un plan (voir ci-dessus). Un autre fait classique est qu'au
voisinage d'un point où la courbure de Gauss de $f$ ne s'annule pas, il y a un espace de dimension infinie
de déformations isométriques infinitésimales, correspondant aux solutions locales d'une EDP linéaire d'ordre $2$
de type elliptique ou hyperbolique.
\medskip

Compte tenu de la nature métrique de la définition, le résultat suivant est un peu surprenant
\begin{thm} (Darboux-Sauer) La rigidité infinitésimale est projectivement invariante. 
Si $f:S\to\R^3\subset\R P^3$ est une immersion, et $\alpha\in\mathrm{PGL}_4(\R)=\mathrm{Aut}(\R P^3)$ 
est une transformation projective telle que $\alpha\circ f(S)\subset \R^3$, alors $\alpha\circ f$ 
est infinitésimalement rigide si et seulement si $f$ l'est.
\end{thm}
\dem
La preuve consiste à remarquer que le groupe $\mathrm{PGL}_4(\R)$ agit sur $\R P^7=P(\R^4\oplus\R^4)$ 
en respectant la quadrique $Q$ et ses sous-espaces projectifs $P_1$, $P_2$, et en commutant à la projection
$\pi_1:Q\setminus P_2\to P_1$~: 
il suffit de faire agir $A\in \mathrm{GL}_4(\R)$
par $\mathrm{diag}(A,{}^tA^{-1})\in SO(q)$ sur $\R^4\oplus\R^4$, et de
noter que cette action préserve $Q_+$, 
car elle en fixe l'élément $(0,0)_+=P(\R^4\oplus 0)=P_1$.\cqfd

\begin{remq}
Cette preuve montre que les ensembles de déformations isométriques infinitésimales
d'une immersion $f$ et de sa transformée projective $\alpha\circ f$ se correspondent bijectivement, 
et que cette correspondance est linéaire (ce sont des espaces vectoriels). Cela tient au fait que 
l'action de $\mathrm{PGL}_4(\R)$ ci-dessus respecte la structure de fibré vectoriel 
sur $\pi_1:Q\setminus P_2\to P_1$~: cette action  s'identifie à celle de $\mathrm{PGL}_4(\R)$ 
sur le fibré cotangent $T^* \R P^3 \to \R P^3$ (isomorphe au fibré tangent).
\end{remq}

\begin{remq}
Darboux démontre ce résultat dans \cite[Ch. IV, n° 900]{Darboux}, et il prouve aussi l'invariance par dualité projective
(ou plutôt ``transformation par polaires réciproques''). Sauer reprend l'étude (trente ans plus tard) dans \cite{Sauer}
en utilisant les coordonnées (projectives) de Plücker sur l'espace des droites affines de $\R^3$. Voir \cite{Izmestiev}, 
pour une exposition récente. L'introduction de la quadrique $Q$ simplifie quelque peu les calculs, 
et fait apparaître le groupe de symétries $\mathrm{Aut}(Q)\simeq PGO(4,4)$ (de dimension 28), 
déjà présent en filigrane dans \cite[ch.IV, n° 903, 907]{Darboux}. 
On peut aussi appliquer cette méthode de preuve dans le cas discret 
(polyèdres et ``armatures''\footnote{pour ``frameworks''}), 
ainsi qu'en plus grande dimension. 
\end{remq}


\section{Triplets de Darboux}\label{triplets}

Soit $f:S\to\R^3$ une immersion, et $g:S\to\R^3$ une déformation isométrique infinitésimale (toutes deux lisses). 
L'équation $df\cdot dg=0$ entraîne qu'il existe une unique application $h:S\to\R^3$ telle que 
$$dg=h\times df,$$
car une rotation de $\R^3$, même infinitésimale, est déterminée par sa restriction à un plan
vectoriel quelconque. 
\begin{defn} On appellera {\em triplets de Darboux} sur une surface $S$ les triplets d'applications lisses 
$(f,g,h)$ de $S$ vers $\R^3$ vérifiant cette relation.
\end{defn}

\def\tg{\tilde{g}}
Si on pose avec Darboux $\tg=g-h\times f$, on obtient $d\tg=f\times dh$, i.e. $(h,\tg,f)=A(f,g,h)$ est un triplet de Darboux.
En particulier $\tg:S\to\R^3$ est une déformation isométrique infinitésimale de $h:S\to\R^3$, 
même si cette dernière n'est pas nécessairement une immersion. Ainsi, un déplacement infinitésimal $g=a\times f+b$
($a,b\in\R^3$), donne les applications constantes $h=a$ et $\tg=b$.

De même, la symétrie de l'équation $df\cdot dg=0$ entraîne que, au moins dans l'ouvert $S'$ de $S$ où $dg$ est de rang $2$, 
il existe une unique application $h^*:S'\to\R^3$ telle que $df=h^*\times dg$, et $(g,f,h^*)=D(f,g,h)$ est un triplet de Darboux sur $S'$.

Il est clair que $A$ est une involution, et que $D$ en est une
sur les triplets de Darboux $(f,g,h)$ où $f$ et $g$ sont des immersions, ce qui équivaut à $h\not\in\im df$. 
On va maintenant préciser la nature de $h^*$.

Tout d'abord, le fait que $h\times df=dg$ est fermée (comme $1$-forme à valeurs vectorielles)
entraîne que $\im dh$ est contenue dans $\im df$. En effet, dans des coordonnées locales sur $S$
on a $\partial_1 h\times\partial_2 f=\partial_2 h\times\partial_1 f$, et en décomposant chaque $\partial_i h$
dans la base $\partial_1 f,\partial_2 f,\nu_f$, où $\nu_f$ engendre $(\im df)^\bot$, on obtient aisément $\partial_i h\in\im df$.

Des égalités $$df=h^*\times dg=h^*\times(h\times df)=(h^*\cdot df)h-(h^*\cdot h)df,$$
il découle (en appliquant $h\times$) que $$dg=-(h^*\cdot h)dg,$$
d'où $$h^*\cdot h=-1,\ \ h^*\cdot df=0.$$
Comme $\im dh$ est contenue dans $\im df$, on a $h^*\cdot dh=0$, 
de sorte que si $h$ est une immersion, $h^*$ est la surface polaire de $h$ par rapport à la quadrique
$x\cdot x=-1$. Autrement dit, le plan affine tangent à $h$ en $p\in S$ a pour équation $h^*(p)\cdot x+1=0$,
et en particulier $h^*(p)$ n'est pas défini si ce plan contient $0\in\R^3$ 
($h^*$ doit alors être interprété projectivement).

Inversement on peut étendre l'involution $D$ à cette situation 
\begin{lemm} Soit $(f,g,h)$ un triplet de Darboux tel que $h$
soit une immersion dont les plans affines tangents 
évitent l'origine $0$. Alors $(g,f,h^*)$ est un triplet de Darboux, 
que l'on prendra pour définition de $D(f,g,h)$.
\end{lemm}
\dem Le fait que $h\times df=dg$ est fermée entraîne maintenant $\im df\subset\im dh$, donc $h^*\cdot df=0$.
Dès lors $h^*\times dg=h^*\times(h\times df)=(h^*\cdot df)h-(h^*\cdot h)df = df$.
Cette définition est consistante, i.e. les deux définitions de $D(f,g,h)$ coïncident quand elles s'appliquent
toutes deux, comme il résulte des calculs précédant le lemme.\cqfd

\def\th{\tilde{h}}
On peut donc appliquer l'involution $D$ au triplet $A(f,g,h)=(h,\tg,f)$ dès lors que 
$f$ est une immersion à plan affine tangent évitant l'origine, 
obtenant le triplet $DA(f,g,h)=(\tg,h,f^*)$, et enfin on peut considérer 
$ADA(f,g,h)=(f^*,\th,\tg)=(f^*,h-f^*\times\tg,\tg)$ 
qui interviendra dans la section suivante.

\begin{remq} Dans l'étude locale des déformations infinitésimales isométriques $g$ d'une immersion $f:S\to\R^3$,
l'hypothèse que $g$ est une immersion n'est pas réellement restrictive, car on peut modifier trivialement
$g$ en lui ajoutant une rotation infinitésimale $a\times f$, et pour un $a\in\R^3$ générique $g+a\times f$
sera une immersion au voisinage de tout point $p\in S$ fixé à l'avance (il suffit que $h+a\not\in\im df$).
De même, l'hypothèse que les hyperplans affines tangents à $f$ évitent l'origine peut être satisfaite
au voisinage de $p$ en remplaçant $f$ par $f+b$ pour un $b\in\R^3$ générique.
\end{remq}

\section{Calcul de $\phi_\pm$}  \label{calcul}

\def\hphi{\hat{\phi}}
Soit $\phi=(f:1:g:-f\cdot g):S\to Q$ l'immersion totalement isotrope associée à une immersion $f:S\to\R^3$
et à sa déformation isométrique infinitésimale $g$, et $h:S\to\R^3$ vérifiant $dg=h\times df$.  
Alors le plan projectif tangent à $\phi$ en $p\in S$ est $\tau_p\phi=P(V_p)$, où 
$V_p=\R\psi(p)\oplus\im d\psi(p)$ est totalement
isotrope de dimension $3$ dans $\R^{4,4}$ et
$$\psi=(f,1,g,-f\cdot g) :S\to\R^{4,4}.$$
On vérifie alors sans peine à partir de $dg=h\times df$ et $\tg=g-h\times f$ 
que $V_p^\bot$ a pour équations en $(x,s,y,t)\in\R^{4,4}$ (et en sous-entendant l'évaluation en $p$)
$$ (y-h\times x-\tg s)\cdot df=0, \ \ 
 (x-sf)\cdot g +y\cdot f +t =0,$$
que l'on peut réécrire
$$ y =h\times x+\tg s+\lambda\nu_f ,\ \ 
   t =-\tg\cdot x -\lambda\nu_f\cdot f,$$
où $\nu_f\in\R^3$ engendre $(\im df)^\bot$ et $\lambda\in\R$.
Si $(x,s,y,t)\in V_p^\bot$, on a donc $x\cdot y+st=\lambda\nu_f\cdot(x-sf)$,
de sorte que les deux sous-espaces totalement isotropes maximaux de $\R^{4,4}$ contenant $V_p$
sont donnés par les équations additionnelles $\lambda=0$ et $\nu_f\cdot(x-sf)=0$, soit encore
$$ V^+_p \left\{ \begin{array}{lll} y & = & h\times x+\tg s \\ t & = & -\tg\cdot x \end{array} \right.$$
$$ V^-_p \left\{ \begin{array}{lll} 
y & = & h\times x+\tg s +\lambda\nu_f \\
t & = & -\tg\cdot x -\lambda\nu_f\cdot f \\
0 & = & \nu_f\cdot(x-sf)
\end{array} \right.$$
En particulier, avec les paramétrages $(.,.)_\pm$ d'ouverts affines de $Q_\pm$ 
par $\R^3\oplus\R^3$ donnés section \ref{dii} on a $$\phi_+=(h,\tg)_+\ .$$

Pour donner à $V_p^-$ une allure plus présentable, commençons par supposer que $\nu_f\cdot f\neq 0$ en $p$.
On peut alors définir $f^*=-\nu_f/(\nu_f\cdot f)$, polaire de $f$ par rapport à la quadrique $x\cdot x+1=0$ (voir la section précédente), 
et les équations définissant $V^-_p$ peuvent s'écrire
$$V^-_p \left\{ \begin{array}{lll} 
y & = & h\times x+\tg s -\lambda(\nu_f\cdot f) f^* \\
t & = & -\tg\cdot x -\lambda\nu_f\cdot f \\
s & = & -f^*\cdot x
\end{array} \right.$$
ou encore, en éliminant $\lambda\nu_f\cdot f$ entre les deux premières équations
$$V^-_p \left\{ \begin{array}{lll} 
y & = & h\times x+\tg s +(\tg\cdot x+t) f^* \\
s & = & -f^*\cdot x
\end{array} \right.$$
qui peut aussi s'écrire
$$V^-_p \left\{ \begin{array}{lll} 
y & = & (h-f^*\times\tg)\times x + f^* t \\
s & = & -f^*\cdot x
\end{array} \right.$$
Noter que $(\tg,h,f^*)$ est un triplet de Darboux d'après la section précédente, donc aussi 
son image $A(\tg,h,f^*)=(f^*,\th,\tg)=(f^*,h-f^*\times\tg,\tg)$.
En particulier, on a $$\phi_-=(\th,f^*)_-$$ là où $f^*$ est définie.

Quand $\nu_f\cdot f$ s'annule en $p$, les équations définissant $V_p^\bot$ données plus haut deviennent
$$V^\bot_p \left\{ \begin{array}{lll} 
y & = & h\times x + s\tg +\lambda\nu_f \\
t & = & -\tg\cdot x
\end{array} \right.$$
et la forme quadratique $x\cdot y+st$ vaut $\lambda\nu_f\cdot x$ sur ce sous-espace,
de sorte que $V_p^-$ est donné par l'équation additionnelle $$\nu_f\cdot x=0.$$

Remarquons que $\nu_f\cdot f$ s'annule identiquement si et seulement si l'immersion $f$ 
prend ses valeurs dans un cône de sommet $0$. Cependant, même dans ce cas le sous-espace $V^-_p$ précédent
est parfaitement défini, et son image projective est $\phi_-(p)$.


\section{Matrices de Zorn \cite{vdBS}} \label{octo}

\def\ZZ{\mathcal{Z}}
On considère l'algèbre $\ZZ$ des matrices de Zorn sur $\R$
$$\begin{pmatrix} a & x \\ y & b \end{pmatrix}, \ \ a,b\in\R, x,y\in\R^3$$ 
de produit
$$\begin{pmatrix} a & x \\ y & b \end{pmatrix}\begin{pmatrix} a' & x' \\ y' & b' \end{pmatrix}=
\begin{pmatrix} aa'+x\cdot y' & ax'+b'x-y\times y' \\ a'y+by'+x\times x' & x'\cdot y+bb' \end{pmatrix}$$
Alors le ``déterminant'' $$N(\begin{pmatrix} a & x \\ y & b \end{pmatrix})=ab-x.y$$ est multiplicatif, car son
défaut de multiplicativité avec le produit ``standard'', i.e. sans les produits vectoriels, est 
$(x\cdot y')(x'\cdot y)-(x\cdot y)(x'\cdot y')$, qui est l'opposé du seul terme supplémentaire 
$(x\times x')\cdot(y\times y')$ introduit par ces produits vectoriels.

Comme $\ZZ$ a clairement une unité, il en résulte que $\ZZ$ est une algèbre d'octonions déployée (sur $\R$).
En particulier les sous-espaces totalement isotropes maximaux de $\ZZ$ sont les annulateurs gauches
et droits des éléments isotropes non nuls de $\ZZ$ (voir section \ref{trial}).

Par exemple, on vérifie aisément que l'annulateur à gauche de 
$$\begin{pmatrix} 1 & u \\ v & u\cdot v \end{pmatrix} \in\ZZ$$
n'est autre que l'image du sous-espace de $\R^{4,4}$ noté $(u,v)_-$ section \ref{dii}
par l'application 
$$\begin{array}{rccc}
\rho: & \R^{4,4} & \to & \ZZ \\
 & (x,s,y,t) & \mapsto & \begin{pmatrix} s & x \\ y & -t \end{pmatrix},
\end{array}$$
qui est une anti-isométrie, i.e. $N(\rho(x,s,y,t))=-x\cdot y -st$.

L'annulateur à droite du même élément de $\ZZ$ est quant à lui l'image par $\rho\sigma$ de $(v,u)_-$ (sic),
avec $\sigma(x,s,y,t)=(y,s,x,t)$. Un calcul sans difficulté montre alors que si $u\cdot v\neq 0$,  
ce sous-espace est l'image par $\rho$ de $(-u/(u\cdot v),-v/(u\cdot v))_+$. Plus généralement $(u,v)_+$ est
l'annulateur à droite de 
$$\begin{pmatrix} u\cdot v & -u \\ -v & 1 \end{pmatrix}=\overline{\begin{pmatrix}  1 & u \\ v & u\cdot v \end{pmatrix}} \in\ZZ,$$
où $Z\mapsto\overline{Z}$ désigne l'anti-automorphime de conjugaison dans $\ZZ$.

En particulier, en identifiant $\R^{4,4}$ à $\ZZ$ par l'application $\rho$, 
on peut paramétrer $Q_+$ par $Q\subset P(\ZZ)$ via l'application $\R Z_+\mapsto \ker L_{Z_+}$,
($Z_+\in\ZZ\setminus\{0\},\ N(Z_+)=0$) 
et de même pour $Q_-$ via $\R Z_-\mapsto \ker R_{Z_-}$, où $L_Z$ et $R_Z$ désignent les multiplications à gauche 
et à droite par $Z$ dans $\ZZ$. 

\begin{lemm} \label{incid} \cite{vdBS} La relation d'incidence $\dim(\ker L_{Z_+}\cap\ker R_{Z_-})=3$ entre $Q_+$ et $Q_-$
équivaut à $Z_-Z_+=0$.
\end{lemm}
\dem Il suffit de le vérifier dans un ouvert dense de $Q_+\times Q_-$, car les deux relations définissent des sous-variétés
compactes de $Q_+\times Q_-$ (toutes deux fibrées en $\R P^3$ sur $Q_\pm$).
Or un calcul rapide montre que $(a,b)_+$ et $(a',b')_-$ sont incidents si et seulement si $$a'-a=b\times b',$$
et d'après ce qui précède ils s'identifient via $\rho$ à $\ker L_Z$ et $\ker R_{Z'}$ où 
$$Z'=\begin{pmatrix} 1 & a' \\ b' & a'\cdot b' \end{pmatrix}, \ \ Z=\begin{pmatrix} a.b & -a \\ -b & 1 \end{pmatrix}.$$
Il n'y a plus qu'à vérifier que la relation ci-dessus équivaut à $Z' Z=0$ dans l'algèbre $\ZZ$.
\cqfd


\section{Trialité différentielle} \label{trialdiff}

Soit $\OO$ une algèbre d'octonions deployés sur $\R$, $Q\subset P(\OO)$ la quadrique des droites isotropes de $\OO$,
$\phi:S\to Q$ une immersion totalement isotrope d'une surface. Alors, quitte à restreindre à un ouvert 
(ou bien passer à un revêtement double) de $S$, $\phi$ se relève en $\psi:S\to\OO\setminus\{0\}$, et pour tout $p\in S$,
le sous-espace $V_p=\R\psi(p)+\im d\psi(p)$ de $\OO$ est totalement isotrope de dimension $3$.

Les deux sous-espaces totalement isotropes maximaux contenant $V_p$ sont $\ker L_{\psi_+(p)}$ et $\ker R_{\psi_-(p)}$
pour des applications lisses $\psi_\pm:S\to\OO\setminus\{0\}$ (là encore quitte à travailler localement), de sorte qu'en termes de
la structure d'algèbre de $\OO$
$$\begin{array}{rllrllrll} 
\psi_+\psi & = & 0, & \psi\psi_- & = & 0, & \ \ \psi_-\psi_+ & = & 0\\
\psi_+ d\psi & = & 0, &  d\psi\,\psi_-& = & 0, &&&
\end{array}$$
d'où aussi (Leibniz)
$$\begin{array}{rllrllrll} 
d\psi_+\,\psi & = & 0, & \psi d\psi_- & = & 0, & d\psi_-\,\psi_+ +\psi_- d\psi_+& = & 0
\end{array}.$$
La trialité différentielle se résume alors à l'énoncé
\begin{prop} \label{prop:trialdiff}
On a en fait $d\psi_-\,\psi_+ = \psi_- d\psi_+ =  0$.
\end{prop}
\dem Si $X,Y,Z$ sont des champs de vecteurs tangents à $S$, $XY\psi$ est orthogonal à $\psi$, car
$V=\R\psi+\im d\psi$ étant totalement isotrope, 
$$0=X\langle Y\psi,\psi\rangle=\langle XY\psi,\psi\rangle+\langle Y\psi,X\psi\rangle=\langle XY\psi,\psi \rangle.$$
De plus $XY\psi$ est orthogonal à $Z\psi$, car en dérivant $\langle Y\psi,Z\psi \rangle= 0$ on obtient
$$\langle XY\psi,Z\psi\rangle+\langle Y\psi,XZ\psi \rangle =0,$$
ce qui montre que la quantité $\langle XY\psi,Z\psi\rangle$ est antisymétrique en $Y,Z$, alors
qu'elle est symétrique en $X,Y$. Par un lemme classique, elle est donc nulle.

En particulier, $XY\psi$ prend ses valeurs dans 
$$V^\bot=(\ker L_{\psi_+}\cap\ker R_{\psi_-})^\bot .$$
Avant de poursuivre la preuve, rappelons le
\begin{lemm} Pour $a\in\OO$ posons $\mathrm{Re}(a)=\langle a,1\rangle$ et $\overline{a}=2\mathrm{Re}(a)1-a$.
Alors $a\mapsto \overline{a}$ est un anti-automorphisme de $\OO$ et $L_a$, $L_{\overline{a}}$ sont adjoints, de même que
$R_a$ et $R_{\overline{a}}$. De plus $L_aL_{\overline{a}}=N(a)\mathrm{Id}=R_{\overline{a}}R_a$. 

Si $a,b\in\OO\setminus\{0\}$ vérifient $ba=0$, on a $a(\OO\overline{b})=\R\overline{b}$ et $(\overline{a}\OO)b=\R\overline{a}$.
\end{lemm}
\dem Voir \cite{vdBS}. Disons simplement qu'elle s'obtient en polarisant 
l'identité multiplicative satisfaite par la norme, ce qui donne pour $a,b,x,y\in\OO$
$$\langle ax,by\rangle + \langle ay,bx\rangle = 2\langle a,b\rangle \langle x,y\rangle,$$
puis en appliquant cette identité un nombre suffisant de fois. Par exemple, on obtient sous les hypothèses du lemme
$a(x\overline{b})=2\mathrm{Re}(ax)\overline{b}$, d'où la seconde assertion 
(l'autre identité s'en déduit par conjugaison).
\cqfd

Revenant à la preuve de la proposition, il résulte du lemme que 
$$V^\bot=\im L_{\overline\psi_+} + \im R_{\overline\psi_-}=\overline{\psi_+}\OO+\OO\overline{\psi_-},$$ 
puis que $\psi_+(XY\psi)\in \R\overline{\psi_-}$ et $(XY\psi)\psi_-\in \R\overline{\psi_+}$, car
$\psi_-\psi_+=0$.

On obtient ainsi deux champs de formes bilinéaires symétriques sur $S$, notés $\II_\mp \in S^2T^*S$, 
bien définis à multiplication près par des fonctions ne s'annulant pas, i.e. deux structures pseudo-conformes 
(pouvant dégénérer) sur $S$.
 
Soit $X$ un champ de vecteurs sur $S$ et montrons que $\psi_- X\psi_+$ s'annule en $p\in S$.
On peut supposer $X\psi_+(p)\neq 0$, et choisir un champ de vecteurs $Y$ sur $S$ tel que $Y(p)\neq 0$
et $\II_-(X(p),Y(p))=0$, i.e. $\psi_+(p)(XY\psi)(p)=0$. Alors en dérivant dans la direction $X(p)$ l'identité
$\psi_+ Y\psi=0$, on obtient $X\psi_+(p) Y\psi(p)=0$, de sorte que $\psi(p)$ et $Y\psi(p)$ sont dans le noyau
de $L_{X\psi_+(p)}$. Mais ils sont aussi dans le noyau de $R_{\psi_-(p)}$ (par définition de $\psi_-$), et comme
l'intersection de ces noyaux est de dimension impaire elle est de dimension $3$. 
Finalement $\psi_-(p) X\psi_+(p)=0$ d'après le lemme \ref{incid}.  \cqfd

Il en résulte que pour tout $p\in S$, tout point du plan projectif tangent $\tau_p\phi_+ \subset Q_+$
(resp. de $\tau_p\phi_-  \subset Q_-$) est incident à $\phi(p)\in Q$ et $\phi_-(p)\in Q_-$ (resp. à $\phi(p)\in Q$ et $\phi_+(p)\in Q_+$),
tout comme (par définition) tout point de $\tau_p\phi \subset Q$ est incident à $\phi_+(p)\in Q_+$ et $\phi_-(p)\in Q_-$. 

En particulier, si $\phi_+$ (resp. $\phi_-$) est une immersion, $\phi_-$ et $\phi$ (resp. $\phi$ et $\phi_+$) 
sont ses deux ``applications de Gauss''. 


\section{Identification des transformations de Darboux}\label{darbouxZorn}

Il est naturel au vu des deux sections qui précèdent d'identifier chacune des quadriques $Q,Q_+,Q_-$ 
avec la quadrique $Q_0\subset P(\ZZ)$ au moyen des isomorphismes
$$\begin{array}{ccc} \rho : Q\to Q_0 & \rho_+ : Q_+ \to Q_0 &  \rho_- : Q_- \to Q_0 \\
(x:s:y:t)\mapsto \R\begin{pmatrix}s&x\\y&-t\end{pmatrix} \ &
\rho^{-1}(\ker L_{Z_+}) \mapsto \R Z_+ \ & \rho^{-1}(\ker R_{Z_-}) \mapsto \R Z_-
\end{array}
$$
(rappelons (section \ref{octo}) que $\ZZ$ désigne l'algèbre des matrices de Zorn sur $\R$).
On obtient ainsi des isomorphismes ``trialitaires''
$$ \theta_\pm =\rho^{-1}\rho_\pm: Q_\pm\to Q.$$
Associons à chaque triplet de Darboux $T=(f,g,h)$ sur $S$ l'application
$$\phi_T=(f:1:g:-f\cdot g):S\to Q.$$

\begin{prop}
Il existe des applications $A$ et $D$ de l'ensemble des immersions totalement isotropes $\phi:S\to Q$ 
vers les applications totalement isotropes de $S$ dans $Q$, telles que
si $T=(f,g,h)$ est un triplet de Darboux pour lequel $\phi_T$ est une immersion 
(i.e. $f$ est une immersion), $A(\phi_T)=\phi_{A(T)}$ et que si de plus $D(T)$ est défini, 
$D(\phi_T)=\phi_{D(T)}$.
\end{prop}
\dem
En prenant $\sigma(x:s:y:t)=(y:s:x:t)$, il est clair que $D(\phi)=\sigma\circ\phi$  convient. 
Rappelons que $D(f,g,h)=(g,f,h^*)$ n'a été défini que si $g$ est une immersion (donc $f$ aussi) ou si $h$ 
est une immersion dont les plans affines tangents évitent l'origine.
Pour traiter le cas de la transformation $A$, démontrons le
\begin{lemm} \label{pmvsAD} Soit $T=(f,g,h)$ un triplet de Darboux, avec $f$ une immersion. Alors
$$\theta_+\circ(\phi_T)_+ = (h:-h\cdot\tg:\tg:1)=c\circ\phi_{A(T)},$$
où $c$ est la transformation projective induite par la conjugaison, i.e. $c(x:s:y:t)=(x:t:y:s)$.
\end{lemm}
\dem En remarquant que $c$ est induite via $\rho$ par la conjugaison $Z\mapsto \bar{Z}$ de $\ZZ$,
la formule peut s'écrire $\rho_+\circ(\phi_T)_+=\overline{\rho\circ\phi_{A(T)}}$, ce qui équivaut à
$(\phi_T)_+=\ker L_{\overline{\rho\circ\phi_{A(T)}}}$. 
Mais cette égalité résulte
du calcul de $\phi_+$ section \ref{calcul} 
et de l'identification de $(a,b)_+$ en termes de $\ZZ$ section \ref{octo}. 
\cqfd

Le lemme montre que l'on peut prendre $A(\phi)=c\circ\theta_+\circ\phi_+$,
et la proposition est démontrée.
\cqfd

Pour alléger les notations, identifions les quadriques $Q_\pm$ à $Q$ 
via les isomorphismes trialitaires $\theta_\pm$, ce qui permet de voir 
$\phi_\pm$ comme applications de $S$ vers $Q$. 
La trialité différentielle démontrée dans la section \ref{trialdiff} précédente se résume
alors à la
\begin{prop}
Si $\phi:S\to Q$ est une immersion totalement isotrope telle que $\phi_+$ soit
une immersion, on a $\phi_{++}=\phi_-$ et $\phi_{+-}=\phi$. De même, si $\phi_-$ est une
immersion on a $\phi_{--}=\phi_+$ et $\phi_{-+}=\phi$.\cqfd
\end{prop}
\begin{cor} Si $\phi:S\to Q$ est une immersion totalement isotrope telle que $\phi_\pm$ soient
aussi des immersions, les opérations $(\ )_\pm$ sur $\{\phi,\phi_+,\phi_-\}$ sont les deux
générateurs d'une action de $\Z/3\Z$. \cqfd
\end{cor}
Le fait que $c$ est induit par un anti-automorphisme de $\ZZ$ entraîne que 
$(c\circ\phi)_+=c\circ\phi_-$, et la proposition entraîne que $A$ est une involution sur
les immersions $\phi$ telles que $A(\phi)$ soit une immersion. 

De même que $c$, l'involution $\sigma$ telle que $D(\phi)=\sigma\circ\phi$ 
est induite via $\rho$ par un anti-automorphisme de $\ZZ$, la transposition.
On alors $D\circ A(\phi)=\sigma\circ c\circ\phi_+$, et comme $\sigma\circ c$ est un automorphime involutif
de $\ZZ$, on voit que $(D\circ A)^2(\phi)=\phi_{++}=\phi_-$ dès que cela a un sens.
En particulier le corollaire précédent entraîne
\begin{cor}[Théorème de Darboux] La transformation $(D\circ A)^6$ est l'identité
sur l'ensemble des immersions totalement isotropes $\phi:S\to Q$ telles que $\phi_\pm$ sont aussi des immersions.\cqfd
\end{cor}
On voit aussi que l'orbite de $\phi$ sous l'action du groupe $D_{12}$ diédral d'ordre $12$ engendré par $A,D$ est
constituée de $\phi,\phi_+,\phi_-$ et de leurs images par les éléments non triviaux $\sigma,c,\sigma\circ c$
du groupe d'ordre $4$ engendré par les involutions projectives $\sigma$ et $c$ de $Q$.
En particulier l'involution centrale $(D\circ A)^3$ correspond à l'involution projective $\sigma\circ c$, 
donnée par $(x:s:y:t)\mapsto (y:t:x:s)$.


\section{Dégénérescences}\label{degen}

Il importe de comprendre les cas de dégénérescence, où l'une des applications tangentes $T_p\phi_\pm$ est de rang $0$ ou $1$
dans un ouvert non vide de points $p\in S$ (dans l'étude locale qui suit, on remplacera $S$ par cet ouvert). 

Avant de commencer, remarquons que pour l'étude locale d'une immersion totalement isotrope $\phi:S\to Q$,
on peut se ramener (quitte à restreindre $S$) par un automorphisme de $Q$ au cas où 
$\phi=(f:1:g:-f\cdot g)$, pour une immersion $f:S\to\R^3$ et sa déformation isométrique infinitésimale $g:S\to\R^3$.

Dans ce cas on a vu section \ref{dii} que les cas où $\phi_+$ ou $\phi_-$ est constante correspondent
respectivement aux déformations triviales (déplacement infinitésimaux) de l'immersion $f:S\to\R^3$ 
et aux ``déformations normales'' de $f$ lorsque son image est contenue dans un plan.

Supposons (la différentielle de) $\phi_+=(h,\tg)_+$ de rang $1$ en tout point de $S$. 
C'est équivalent à $h$ de rang $1$ en tout point de $S$,  car $d\tg=f\times dh$.

On peut alors choisir des coordonnées locales $(x_1,x_2)$ sur $S$ telles que $\partial_1$ engendre $\ker dh$,
d'où $\partial_1 h=\partial_1\tg=0$ et $\partial_1 f\times\partial_2 h=\partial_2f\times\partial_1h=0$.

On en déduit que $h=h(x_2)$ et $\partial_1 f=\lambda(x_1,x_2) h'(x_2)$, pour
une fonction $\lambda$ ne s'annulant pas. Quitte à remplacer $x_1$ 
par $x'_1=\int_0^{x_1} \lambda(t,x_2)\,dt$, on peut supposer que $\lambda=1$, i.e. $\partial_1 f=h'(x_2)$.
On a alors $$f(x_1,x_2)=f_0(x_2)+x_1 h'(x_2),$$
autrement dit $f$ est {\it réglée}. L'intégration de $\partial_1 g=h\times\partial_1 f=h\times h'$ montre alors que
$$g(x_1,x_2)=g_0(x_2)+x_1\, h(x_2)\times h'(x_2)$$
est aussi réglée, de sorte que l'immersion $\phi:S\to Q$ elle-même est réglée, i.e. localement de la forme
$$\phi(x_1,x_2)=\R(\psi_0(x_2)+x_1\psi_1(x_2)),$$
où $\psi_0,\psi_1$ sont telles que $\psi_0(x_2),\psi_1(x_2),\psi'_0(x_2),\psi'_1(x_2)$ 
engendrent pour tout $x_2$ un sous-espace totalement isotrope de $\R^{4,4}$ de dimension au moins $3$,
qui coïncide avec $\phi_+(x_2)$ quand il est de dimension $4$.

Inversement, ces équations définissent des déformations isométriques infinitésimales d'immersions réglées, 
pour lesquelles l'image de $\phi_+$ est une courbe, dès lors que $g'_0=h\times f'_0$
et que $f$ est une bien une immersion, i.e. $\partial_1 f\times\partial_2 f=h'\times (f'_0+x_1 h'')$
ne s'annule pas dans le domaine des $(x_1,x_2)$ considérés. 
\medskip

Envisageons maintenant le cas où $\phi_-$ est de rang $1$ en tout point de $S$.
On peut supposer comme précédemment (travaillant localement) que $\phi=(f:1:g:-f\cdot g)$, 
pour une immersion $f:S\to\R^3$ et sa déformation isométrique infinitésimale $g:S\to\R^3$.
On peut aussi supposer que les plans affines tangents de $f$ évitent l'origine, 
quitte à appliquer un nouvel automorphisme de $Q$ relevant une translation de $\R^3$.
On a alors d'après la section \ref{calcul} $\phi_-=(\th,f^*)_-$.

Comme $d\th=\tg\times df^*$, l'hypothèse équivaut à $f^*$ de rang $1$
(rappelons que $A\circ D\circ A(f,g,h)=(f^*,\th,\tg)$ est un triplet de Darboux (section \ref{triplets})).
Cela signifie que $f$ est une surface {\em développable}, enveloppe d'une famille
à un paramètre de plans.

On peut choisir des coordonnées locales $(x_1,x_2)$ telles que $\partial_1$ engendre le noyau de $df^*$.
Alors, en abusant un peu des notations, $f^*=f^*(x_2)$, et comme $f^*\cdot f=-1$ et $f^*\cdot df=0$, on a $f\cdot df^*=0$ et
$$f(x_1,x_2)=f_0(x_2)+\mu(x_1,x_2) f^*(x_2)\times {f^*}'(x_2),$$
pour une fonction $\mu$ vérifiant $\partial_1\mu\neq 0$ (que l'on peut prendre égale à $x_1$ si on le souhaite),
et une application $f_0$ telle que $f^*\cdot f_0=-1$ et $f^*\cdot f_0'=0$.

L'égalité $d\th=\tg\times df^*$ entraîne d'abord que $\partial_1\th=0$, puis que $\partial_1\tg\times \partial_2 f^*=0$,
de sorte que 
$$\tg(x_1,x_2)=\tg_0(x_2)+\lambda(x_1,x_2) {f^*}'(x_2)$$
pour une fonction réelle $\lambda$.

En  se rappelant que $\th=h-f^*\times\tg$, l'égalité $h=\th+f^*\times\tg$ s'écrit
$$h(x_1,x_2)=h_0(x_2)+\lambda(x_1,x_2) f^*(x_2)\times {f^*}'(x_2)$$
pour $h_0=\th + f^*\times\tg_0$, que l'on reporte dans $g=\tg+h\times f$
(en omettant les arguments pour plus de lisibilité), obtenant
$$g=\tg_0 +\lambda\, {f^*}' + (h_0+\lambda\, f^*\times {f^*}')\times(f_0+\mu\, f^*\times {f^*}'),$$
qui se réduit à
$$g(x_1,x_2)=g_0(x_2) + \mu(x_1,x_2) h_0(x_2) \times (f^*(x_2)\times {f^*}'(x_2))$$
en utilisant $f^*\cdot f_0=-1$ et ${f^*}'\cdot f_0$.

Quant à la condition de fermeture $\partial_1 h\times \partial_2 f=\partial_2 h\times \partial_1 f$,
elle se réduit après un calcul fastidieux à une seule équation (les deux membres sont proportionnels à $f^*$)
$$a\, \partial_1(\lambda\mu) = b\, \partial_1\lambda + c\, \partial_1\mu,$$
pour des fonctions $a,b,c$ de $x_2$ données par
$$a=\det(f^*,{f^*}',{f^*}''),\ b=f'_0\cdot {f^*}',\ c=h'_0\cdot {f^*}' \ .$$
Ceci équivaut clairement à l'existence d'une fonction $d$ de $x_2$ telle que
$a\lambda\mu = b\lambda + c\mu +d$, soit encore $$\lambda= \frac{c\mu+d}{a\mu-b}.$$
Ainsi à $x_2$ fixé, le vecteur $h=h_0+\lambda\, f^*\times {f^*}'$ s'obtient par application à 
$f=f_0+\mu\,f^*\times {f^*}'$ d'une transformation projective entre les droites qu'ils parcourent.

Des formules donnant $f$ et $g$ il résulte qu'ici encore, l'immersion $\phi$ est réglée, i.e. localement de
la forme $\phi(x_1,x_2)=\R(\psi_0(x_2)+x_1\psi_1(x_2)),$ où $\psi_0,\psi_1$ 
sont telles que $\psi_0(x_2),\psi_1(x_2),\psi'_0(x_2),\psi'_1(x_2)$ 
engendrent pour tout $x_2$ un sous-espace totalement isotrope de $\R^{4,4}$ de dimension au moins $3$, 
qui coïncide avec $\phi_-(x_2)$ quand il est de dimension $4$.
\medskip

De ces calculs on déduit la
\begin{prop}
Supposons $S$ compacte sans bord et $S$, $f:S\to\R^3$, $g:S\to\R^3$ analytiques, 
avec $g$ une déformation isométrique infinitésimale de l'immersion $f$.

Alors l'immersion totalement isotrope $\phi=(f:1:g:-f\cdot g):S\to Q$ est telle que
$\phi_-$ est de rang $2$ (i.e. une immersion) sur un ouvert dense de $S$, et il en va 
de même pour $\phi_+$ si elle est non constante. 
\end{prop}
\dem En effet $f$ n'ayant pas d'ouvert plat, ni même réglé (par analyticité et compacité), 
$\phi_-$ ne peut être constante ou de rang $1$ sur un ouvert d'après ce qui précède. 
De même, si $\phi_+$ est non constante, elle ne peut être de rang $1$ sur un ouvert, pour la même raison.\cqfd 

\begin{remq}
Cohn-Vossen, Rembs, Reshetnyak puis Trotsenko \cite{Trotsenko} 
ont construit des exemples de surfaces de révolution $f:S^2\to\R^3$
de plus en plus régulières (analytiques réelles pour ce dernier) 
et possédant des déformations isométriques infinitésimales non-triviales.
Cependant, comme le fait remarquer Spivak dans \cite{Spivak} aucun de ces exemples n'est analytique réel,
ni même $C^\infty$~: le champ de vecteurs $g$ a un ordre de différentiabilité fini aux pôles de la surface
de révolution. Ainsi, la question de l'existence d'une déformation isométrique infinitésimale non-triviale
d'une surface fermée dans la catégorie analytique réelle, i.e. essentiellement celle d'une immersion
totalement isotrope analytique réelle non-dégénérée d'une surface fermée dans la quadrique $Q$ reste ouverte.
\end{remq}


\section{Les secondes formes fondamentales}\label{secff3}

Si $\phi:S\mapsto Q$ est une immersion totalement isotrope, on peut la relever localement en 
$\psi:S\to\R^{4,4}\setminus\{0\}$, et alors $V=\R\psi+\im d\psi$ est totalement isotrope de dimension $3$.

On a vu section \ref{trialdiff} que pour tous champs de vecteurs $X,Y,Z$ tangents à $S$,
on a $XY\psi \in V^\bot$, d'où $\langle XY\psi ,Z\psi\rangle=0$ et $\langle XY\psi,\psi\rangle=0$, 
ce qui implique que pour tous $X,Y,Z$ on a $\langle XYZ\psi,\psi\rangle=0$. Si $W$ est un quatrième champ de vecteurs tangents à $S$,
on en déduit que $\langle XYZ\psi,W\psi\rangle = -\langle YZ\psi, XW\psi\rangle$ est symétrique en $(W,X,Y,Z)$,
puisque c'est aussi $-\langle WXYZ\psi,\psi\rangle$. 

Mais comme $V^\bot/V\simeq\R^{1,1}$ a deux droites istotropes $\R v_\pm$ , on peut poser
$$XY\psi\mod V = \II_+(X,Y) v_+ + \II_-(X,Y) v_-$$ 
pour deux formes bilinéaires symétriques $\II_\pm \in S^2T^*S$ bien définies à un facteur près.
On peut supposer $\langle v_+,v_-\rangle =1$, et aussi que $V+\R v_\pm$ appartient à $Q_\pm$, pour fixer les signes.

Alors la totale symétrie de $\langle WX\psi, YZ\psi\rangle$ en $(W,X,Y,Z)$ équivaut à la condition que pour tous $X,Y$ on a
$$ \II_+(X,X)\II_-(Y,Y) +\II_+(Y,Y)\II_-(X,X) = 2 \II_+(X,Y)\II_-(X,Y),$$
autrement dit $\II_+$ et $\II_-$ sont orthogonales pour la forme quadratique discriminant ("apolaires")
$$ \mathrm{discr}(B)\in (\wedge^2T^*S)^{\otimes 2} : X\wedge Y \mapsto B(X,X)B(Y,Y)-B(X,Y)^2.$$

Les formes $\II_\pm$ sont à un facteur non nul près les mêmes que celles définies section \ref{trialdiff}
au moyen des octonions déployés $\OO\simeq\R^{4,4}$ par 
$$\psi_+(XY\psi)=\II_-(X,Y)\overline{\psi_-},\ \ \ (XY\psi)\psi_-=\II_+(X,Y)\overline{\psi_+}.$$
On les appellera {\em secondes formes fondamentales} de l'immersion totalement isotrope $\phi$ 
\footnote{La seconde forme fondamentale de l'immersion projective $\phi:S\to Q\to\R P^7$ est une application linéaire
$\II_\phi:N^*_\phi\to S^2T^*S$, et si $\phi$ est une immersion totalement isotrope, l'image de $\II_\phi$ est engendrée
par $\II_+$ et $\II_-$}.

Noter que l'on a aussi 
$$(XY\psi_+)\psi=\II_-(X,Y)\overline{\psi_-},\ \ \ \psi(XY\psi_-)=\II_+(X,Y)\overline{\psi_+},$$
et la trialité différentielle (proposition \ref{prop:trialdiff}) entraîne 
que $$(XY\psi_-)\psi_+=\psi_-(XY\psi_+)=-(X\psi_-)(Y\psi_+)=\II(X,Y)\overline{\psi}$$ 
pour une troisième forme bilinéaire symétrique $\II$ sur $TS$, définie à un facteur 
près et orthogonale aux deux formes $\II_\pm$ précédentes.

En particulier
\begin{prop} Si les formes $\II_+$, $\II_-$, $\II$ sont non nulles en un point $p\in S$, alors 
\begin{itemize}
\item[-] soit elles sont toutes non dégénérées, l'une d'elles est définie,
les deux autres sont indéfinies à directions isotropes mutuellement orthogonales et aussi orthogonales pour la forme définie,
\item[-] soit l'une des formes est de rang $1$, disons $\pm\ell^2$ pour une forme linéaire non nulle $\ell\in T^*S$,
les deux autres sont divisibles par $\ell$, et l'une d'elles au moins est proportionnelle à $\ell^2$.
\end{itemize}
\end{prop}
\dem Il s'agit de comprendre les positions de trois points de $\R P^2=P(\R^{1,2})\simeq P(S^2T^*S)$
qui sont deux à deux orthogonaux vis-à-vis de la conique lisse $C$ d'équation $\mathrm{discr}=0$. 
Si aucun n'est sur la conique, un des points
est intérieur à $C$, les deux autres sur sa droite polaire, donc extérieurs à $C$ , et orthogonaux entre eux.
Or l'orthogonalité de deux formes quadratiques non dégénérées en dimension $2$ équivaut à ce que les directions isotropes 
(éventuellement complexes) de l'une sont orthogonales pour l'autre.

Enfin l'orthogonal (polaire) d'un point $q=[\ell^2]$ de $C$ est la droite tangente $d$ à $C$ en ce point, projectif du sous-espace 
de dimension $2$ des formes quadratiques divisibles par $\ell$. Deux points de $d$ sont orthogonaux si et seulement si l'un d'eux est égal à $q$.  
\cqfd

\begin{lemm} \label{lemmker}
En tout point $p\in S$, $\ker (\II_\pm)_p=\ker T_p\phi_\mp$, et chaque noyau est contenu dans $\ker\II_p$.
\end{lemm}
\dem Prenons comme ci-dessus des relèvements $\psi,\psi_\pm:S\to\OO\setminus\{0\}$ de $\phi,\phi_\pm$ au voisinage de $p$. 
Le sous-espace $V=\R\psi(p)+\im d\psi(p)$ est totalement isotrope de dimension $3$, et $X\in \ker(\II_+)_p$ équivaut à
$(Y\psi)(X\psi_-)=0$ pour tout $Y\in T_pS$, donc à $V(X\psi_-)=0$. Comme l'intersection de deux annulateurs gauches dans $\OO$
est de dimension paire, ceci équivaut à $X\psi_-\in\R\psi_-(p)$, autrement dit à $X\in \ker T_p\phi_-$.
De plus, $X\psi_-\in\R\psi_-(p)$ entraîne clairement $(X\psi_-)(Y\psi_+)=0$ pour tout $Y\in T_pS$, i.e. $X\in \ker\II_p$.

Le cas de $\II_-$ est similaire (ou s'en déduit par conjugaison).
\cqfd

On en déduit la caractérisation suivante des situations où l'une des formes $\II_+$, $\II_-$, $\II$ s'annule.
\begin{prop} Pour une immersion totalement isotrope $\phi:S\to Q$, la forme $\II_+$ (resp. $\II_-)$ 
s'annule au point $p\in S$ si et seulement si $T_p\phi_-=0$ (resp. $T_p\phi_+=0$). On a alors aussi $\II_p=0$.

La forme $\II$ s'annule en $p$ si et seulement si $T_p\phi_-=0$, ou $T_p\phi_+=0$, 
ou bien $T_p\phi_\pm$ sont toutes deux de rang $1$.
Dans ce dernier cas, $(\II_\pm)_p$ sont de rang $1$ et ont même noyau $\R X=\ker T_p\phi_\mp$, 
et les droites projectives tangentes $\tau_p\phi_\pm\subset Q_\pm$
découpent par incidence dans $Q$ la même droite projective, 
à savoir celle ayant pour vecteur tangent $X\phi\in T_pQ$ en $\phi(p)$.
\end{prop}
\dem La première partie de la proposition résulte aussitôt du lemme précédent.

La forme $\II_+$ est nulle en $p\in S$ si et seulement si le produit $(X\psi)_p(Y\psi_-)_p$ (dans $\OO$) est nul
pour tous champs de vecteurs $X,Y$ sur $S$. Mais cela entraîne que les sous-espaces vectoriels 
$V=\R\psi(p)+\im d\psi(p)$ et $V_-=\R\psi_-(p)+\im d\psi_-(p)$ de $\OO$ s'annulent mutuellement, i.e. $VV_-=0$.
Comme l'intersection de deux annulateurs gauches (resp. droits) dans $\OO$ est de dimension paire 
et que $V$ est de dimension $3$ (car $\phi$ est une immersion), ceci force $V_-$ à être de dimension $1$,
i.e. $T_p\phi_-(p)=0$. La réciproque est claire. Le cas de $\II_-$ se traite semblablement (ou s'en déduit par conjugaison).

L'annulation de $\II_p$ équivaut à $V_-V_+=0$ dans $\OO$, où $V_-=\R\psi_-(p)+\im d\psi_-(p)$, $V_+=\R\psi_+(p)+\im d\psi_+(p)$,
et on a soit $\dim V_+=1$ ou $\dim V_-=1$, soit $V_-$, $V_+$ sont de dimension $2$, i.e. $T_p\phi_\pm$ sont de rang $1$,
et $V_-V_+=0$. 

Inversement, supposons $T_p\phi_\pm$ de rang $1$. Les formes $(\II_\pm)_p$ sont aussi de rang $1$ d'après le lemme,
et comme elles sont orthogonales pour le discriminant, elles ont même noyau $\R X\subset T_pS$, 
qui est aussi le noyau commun de $T_p\phi_\pm$. 

Soit $Y\in T_pS\setminus \R X$.  Alors $(X\psi)(Y\psi_-)=(Y\psi)(X\psi_-)=0$, car $X\psi_-\in\R\psi_-(p)$, et
de même $(Y\psi_+)(X\psi)=0$, de sorte que 
$$\ker L_{Y\psi_+}\cap \ker R_{Y\psi_-}\supset \R \psi(p)\oplus \R X\psi$$
est de dimension au moins $2$, impaire, donc $3$, et $(Y\psi_-)(Y\psi_+)=0$ en vertu du lemme \ref{incid}.
Mais ceci signifie que $\II(Y,Y)=0$, et comme $X\in\ker\II$ d'après le lemme \ref{lemmker}, on a bien $\II_p=0$.

Quant à la dernière assertion, elle résulte de ce que le plan vectoriel $W=\R\psi(p)+\R X\psi$ coïncide
pour raison de dimension avec $\ker L_{\psi_+(p)}\cap\ker L_{Y\psi_+}$ et avec $\ker R_{\psi_-(p)}\cap\ker R_{Y\psi_-}$.
Géométriquement, cela signifie que la droite $P(W)\subset Q$ est l'ensemble de points de $Q$ incidents
à tous les points de $\tau_p\phi_+=P(V_+)$ (resp. de $\tau_p\phi_-=P(V_-)$). Comme $P(W)$ est la droite
projective de vecteur tangent $X\phi$ en $\phi(p)$, la proposition est démontrée.
\cqfd

Les possibilités obtenues pour le triplet $(\rg\II_+, \rg\II_-,\rg\II)$ sont résumées dans le tableau
$$\begin{array}{ccc}
\II_+ & \II_- & \II \\
2 & 2 & 2 \\
2 & 1 & 1 \\
1 & 2 & 1 \\
1 & 1 & 0 \\
0 & \star & 0 \\
\star & 0 & 0
\end{array}$$

Supposons maintenant que $\phi=(f:1:g:-f\cdot g)$ 
pour une immersion $f:S\to\R^3$ et sa déformation isométrique infinitésimale $g:S\to\R^3$.
On va déterminer les trois formes bilinéaires (à un facteur près) en termes des applications $f$, $g$ et $h$,
où $dg=h\times df$. Pour cela, on peut supposer sans perte de généralité, quitte à travailler localement sur $S$ et remplacer
$(f,g,h)$ par $(f+v,g+a\times f+b,h+a)$ pour $v,a,b\in\R^3$ génériques, 
que $g$ est une immersion et que les plans affines tangents à $f$ et $g$ évitent l'origine.
On peut alors définir les polaires $f^*$, $g^*$ par $f^*\cdot df=0,\ f^*\cdot f=-1$, 
$g^*\cdot dg=0,\ g^*\cdot g=-1$, et on prend
$$\II_f=f^*\cdot D^2f,\ \ \ \II_g=g^*\cdot D^2g$$ 
pour représentantes des secondes formes fondamentales de $f$ et $g$ à un facteur près.
Noter que $g^*=-h/(h\cdot g)$, puisque $dg=h\times df$ entraîne $\im dg=(\R h)^\bot$.

D'après la section \ref{triplets}, le triplet de Darboux $D(f,g,h)=(g,f,h^*)$ est donné par
$h^*\cdot df=0$, $h^*\cdot h=-1$. On a alors $f^*=-(h\cdot f^*)h^*$, donc 
$(h^*\cdot f)(h\cdot f^*)=1$, et on pose 
$$\II_h=h^*\cdot D^2h.$$
Lorsque $h$ est une immersion, ses plans affines tangents évitent l'origine, 
et $h^*$ coïncide avec la polaire de $h$.

On se place dans le modèle des matrices de Zorn $\ZZ$ pour les octonions déployés.
En utilisant les isomorphismes trialitaires de la section \ref{darbouxZorn}
et les résultats des sections \ref{calcul} et \ref{octo}, on peut identifier
$\psi_+$,$\psi$,$\psi_-$ aux matrices de Zorn
$$\psi_+=\begin{pmatrix} h\cdot\tg & -h \\ -\tg & 1 \end{pmatrix},\ 
\ \psi=\begin{pmatrix} 1 & f \\ g & f\cdot g \end{pmatrix},\ 
\ \psi_-=\begin{pmatrix} 1 & \th \\ f^* & \th\cdot f^* \end{pmatrix}$$
et dès lors les trois formules 
\begin{align*}
-X\psi_+Y\psi & = \II_-(X,Y)\overline{\psi_-} \\
-X\psi Y\psi_- & = \II_+(X,Y)\overline{\psi_+} \\
-X\psi_-Y\psi_+ & =  \II(X,Y)\overline{\psi}
\end{align*}
permettent de vérifier la
\begin{prop}
Dans ces conditions, les formes bilinéaires symétriques $\II_+$, $\II_-$, $\II$ sur $TS$ sont respectivement proportionnelles à 
$\II_f$, $\II_g$, $\II_h$, avec des facteurs non nuls. 
\end{prop}
\dem Il suffit de calculer un coefficient du produit correspondant de matrices de Zorn dans chaque cas, ce qui donne
$\II_+(X,Y)=-Yf\cdot Xf^*=(XYf)\cdot f^*$, donc $\II_+=\II_f$.
De même $\II(X,Y)=Xf^*\cdot Yh=-f^*\cdot(XYh)$, donc $\II=(h\cdot f^*)\II_h$.
Enfin $\II_-(X,Y)=X\tg\cdot Yf=(f\times Xh)\cdot(h^*\times Yg)=(f\cdot h^*)(Xh\cdot Yg)=-(f\cdot h^*)(h\cdot XYg)$,
donc $\II_-=(h\cdot g)(f\cdot h^*)\II_g$. 
\cqfd

On retrouve ainsi certaines observations de Darboux \cite[p. 71-72]{Darboux} qui regroupe notamment les $12$ surfaces
par paquets de $4$ partageant les mêmes ``lignes asymptotiques'', autrement dit à secondes formes fondamentales proportionnelles.
Dans les notations présentes, ce sont $\{f,f^*,\tilde{f},\tilde{f}^*\}$, $\{g,g^*,\tilde{g},\tilde{g}^*\}$ et 
$\{h,h^*,\tilde{h},\tilde{h}^*\}$. 
Du point de vue de l'action du groupe diédral $D_{12}$ engendré par les involutions $A$ et $D$ sur les triplets de Darboux
(voir sections \ref{triplets} et \ref{darbouxZorn}), les douze surfaces 
sont les premières composantes des triplets dans l'orbite de $(f,g,h)$, et les sous-ensembles ci-dessus
sont les premières composantes des triplets dans les trois orbites du sous-groupe engendré par les involutions $ADA$, correspondant à 
$(\ )^*$, et $DAD$, correspondant à $\widetilde{(\ )}$. Noter que ces involutions commutent, leur produit $(DA)^3$ étant une involution.


\section{La variété d'incidence et sa distribution de dimension $6$} \label{varincid}

Grâce à la trialité différentielle des sections précédentes, 
on va reformuler les déformations isométriques infinitésimales de surfaces
comme surfaces intégrales d'un champ de plans tangents (une "distribution")
de dimension $6$ sur une variété de dimension $11$.

Soit $Q_0=Q$, $Q_1=Q_-$, $Q_2=Q_+$ (sic), les trois quadriques en trialité. 
Chaque $Q_i$ est une copie de la quadrique des droites isotropes de $P(\OO)$, 
où $\OO$ est une algèbre d'octonions déployée sur $\R$. On définit alors la variété d'incidence
$$\cI=\{ ([x_i])_{i\in\Z/3} \mid x_i x_{i+1} = 0,\; i \in\Z/3\}\subset Q_0\times Q_1\times Q_2.$$
Géométriquement, $\cI$ s'identifie à la variété des drapeaux formés d'un point de $Q$ et d'un plan projectif contenu dans $Q$ et contenant ce point.
C'est une variété de dimension $11$, fibrée sur chaque $Q_i$ avec fibres isomorphes à la variété $P(T^*\R P^3)$ des drapeaux (point,plan) de $\R P^3$.

On note $q_i:\cI\to Q_i$, $i\in\Z/3$ ces fibrations. Alors $q_{i+1}\times q_{i-1}:\cI\to Q_{i+1}\times Q_{i-1}$ a pour image la sous-variété
d'incidence $J_i\subset Q_{i+1}\times Q_{i-1}$, fibrée en $\R P^3$ sur $Q_{i\pm 1}$, et $q_{i+1}\times q_{i-1}:\cI\to J_i$ est une fibration
de fibre $\R P^2$. On considère le fibré tangent vertical $\xi_i\subset T\cI$ de cette fibration, i.e.
$$\xi_i=\ker Tq_{i+1}\cap \ker Tq_{i-1}.$$

\begin{lemm} Les $\xi_i$ sont en somme directe dans $T\cI$.
\end{lemm}
\dem Il est clair que les $\xi_i$, $i\in\Z/3$ sont deux à deux en somme directe, car pour tout $i$, $\xi_i\cap\xi_{i+1}=\bigcap_j \ker Tq_j=0$.
De même, $\xi_i\oplus\xi_{i+1}$ est contenu dans $\ker Tq_{i-1}$, et comme  $\xi_{i-1}\cap\ker Tq_{i-1}=0$, on a bien la conclusion du lemme.
\cqfd

Ceci amène à considérer le sous-fibré de rang $6$ du fibré tangent $T\cI$ (ou "distribution")
 $$\xi=\bigoplus_{i\in\Z/3} \xi_i \subset T\cI\; .$$

\begin{prop} Les immersions totalement isotropes $\phi:S\to Q$ 
sont en bijection avec les surfaces intégrales $\Phi:S\to \cI$ de $\xi\subset T\cI$ 
(immersions avec $\im T\Phi\subset\xi$) telles
que $\phi=q_0\circ\Phi$ soit une immersion, via $\phi\mapsto \Phi=(\phi,\phi_-,\phi_+)$.
\end{prop}
\dem Ce n'est qu'une reformulation de la trialité différentielle (proposition \ref{prop:trialdiff}). \cqfd

\begin{remq} On peut montrer que la distribution de rang $6$ $\xi$ est complètement non intégrable. 
Plus précisément, $\xi'=[\xi,\xi]$ est de rang $9$, et $[\xi,\xi']=T\cI$ de rang $11$.
En effet, la restriction de $\xi$ aux fibres de $q_i$, i.e. $\xi\cap\ker Tq_i$, 
n'est autre que la structure de contact canonique sur chaque fibre $q_i^{-1}([x_i])\simeq P(T^*\R P^3)$. 
En particulier, $[\xi,\xi]=\xi'$ contient la somme des $\ker Tq_i$, 
qui est de rang $9$ comme on le vérifie aisément, et coïncide en fait avec elle. 
On peut vérifier de même que $[\xi,\xi']=T\cI$.
\end{remq}


\end{document}